\documentclass[12pt,reqno]{amsart}
\makeatletter
\renewcommand{\subsection}{%
  \@startsection{subsection}{2}{\z@}%
  {3.25ex\@plus1ex\@minus.2ex}%
  {1.5ex\@plus.2ex}%
  {\normalfont\bfseries}}
\makeatother
\usepackage{xcolor}
\usepackage[a4paper,margin=1in]{geometry}
\usepackage{dsfont}
\usepackage{amsmath,amssymb,amsthm,mathtools}
\usepackage{bm}
\usepackage{hyperref}
\usepackage{microtype}
\numberwithin{equation}{section}
\newtheorem{theorem}{Theorem}[section]
\newtheorem{lemma}[theorem]{Lemma}
\newtheorem{proposition}[theorem]{Proposition}

\newtheorem{remark}[theorem]{Remark}
\newtheorem{definition}[theorem]{Definition}
\newcommand{\Om}{\Omega}
\newcommand{\R}{\mathbb R}
\newcommand{\C}{\mathbb C}
\newcommand{\cG}{\mathcal G}
\newcommand{\cN}{\mathcal N}
\newcommand{\cU}{\mathcal U}
\newcommand{\cM}{\mathcal M}

\newcommand{\cX}{\mathcal X}
\newcommand{\cB}{\mathcal B}
\newcommand{\ip}[2]{\left\langle #1,#2\right\rangle}
\newcommand{\norm}[1]{\left\|#1\right\|}
\newcommand{\jcp}{J. Chem. Phys.}

\title[Stabilization of Gray-Scott steady states]
{Finite-Dimensional Internal Stabilization of Steady States of the Gray-Scott System}

\author{Aman Patel}
\address{Department of Mathematics, Indian Institute of Technology Kanpur, Kanpur 208016, India}
\email{amanpat@iitk.ac.in}

\author{Mrinmay Biswas}
\address{Department of Mathematics, Indian Institute of Technology Kanpur, Kanpur 208016, India}
\email{mbiswas@iitk.ac.in}

\subjclass[2020]{93D15, 93B52, 35K57, 35Q92}
\keywords{Gray-Scott system, reaction-diffusion system, finite-dimensional feedback, stabilization, Riccati equation, Fattorini-Hautus criterion, unique continuation, steady state}

\begin{document}

\begin{abstract}
    In this article, we establish the local exponential stabilization of every steady state of the Gray-Scott system with homogeneous Neumann boundary conditions, other than the trivial state $(0,1)^{\top}$, at any prescribed decay rate by finite-dimensional feedback controls supported in an arbitrary nonempty open subset of the domain and acting only on the first equation of the system. Our approach does not require the steady state to be spatially constant and accommodates non-self-adjoint linearizations with variable coefficients. The proof combines unique continuation, a nonvanishing property of the steady-state coupling, the Fattorini-Hautus condition, and Riccati-based feedback.
\end{abstract}

\maketitle

\section{Introduction}
    Reaction–diffusion systems are mathematical models that describe the spatiotemporal evolution of the concentrations of one or more chemical species, governed by diffusion and nonlinear reaction terms. These models account for local chemical reactions, through which substances are transformed into one another, together with diffusion processes that cause the substances to spread throughout a spatial domain. Reaction–diffusion systems model the dynamics of a broad range of phenomena encountered in various areas of science and engineering, such as physics, chemistry, biology, ecology, and geology. From a mathematical perspective, such systems are typically formulated as systems of semilinear parabolic partial differential equations. In particular, a general two-species reaction–diffusion system takes the form:
    \begin{align*}    
        \begin{cases}
            u_t = D_u  \Delta u + f(u,v),\\
            v_t = D_v \Delta v + g(u,v),
        \end{cases}
    \end{align*}
    where, $u(t,x)$ and $ v(t,x)$ denote the concentrations of two chemical species, $D_u$ and $ D_v$ are positive constants called diffusion coefficients and the functions $f$ and $ g$ represent the local reaction kinetics. Reaction–diffusion systems admit a rich variety of solution behaviors, including the formation of travelling waves and other wave-like phenomena. They also exhibit self-organized Turing patterns, such as spots, hexagons, stripes, and labyrinthine structures, as well as more intricate phenomena including dissipative solitons~\cite{Pea93, LP68, LS76, Mur01, Mur03, JS72}.
\subsection{The Gray-Scott system and its steady states}
    In this work, we focus on the Gray–Scott model~\cite{GS83, GS84, GS85}, a prototypical reaction–diffusion system for two chemical species, combining diffusion with nonlinear autocatalytic reactions and exhibiting rich pattern formation dynamics. The model describes the interaction of two chemical species, denoted by $U$ and $V$, whose concentrations at time $t$ and position $x$ are denoted by $u(t,x)$ and $v(t,x)$, respectively. These species undergo local chemical reactions and simultaneously diffuse through the surrounding medium. As a result, their concentrations evolve over time and vary across the spatial domain. The resulting dynamics are governed by the following coupled system of semilinear partial differential equations (\cite[Eq. 7.20]{Pert15}): 
    \begin{align}\label{Intro_Gray_Scott_System}
		\begin{cases}{}
			  u_t = D_u \Delta u +u^2v-Au ,\\
			  v_t = D_v \Delta v-u^2v+B(1-v) ,  \\
		  \end{cases}
	\end{align}
    which captures the mechanisms governing diffusion, production, and depletion of each species and gives rise to complex spatiotemporal patterns. Here, the constants $D_u>0$ and $D_v>0$ denote the diffusion coefficients, while $A>0$ represents the decay (or removal) rate of species $U$ and $B > 0$ denotes the feed rate of species $V.$

    The first equation governs the spatio-temporal evolution of the concentration $u$ of species $U$ in a reacting and diffusing medium. The diffusion term $D_u\Delta u$ accounts for spatial spreading due to molecular diffusion. The nonlinear coupling term $u^2v$ represents an autocatalytic chemical reaction in which two molecules of species $U$ interact with one molecule of species $V$ to produce additional $U$, according to the reaction scheme $2U + V \rightarrow 3U.$ In addition to diffusion and production through the autocatalytic reaction, the equation contains a linear decay term $-Au$, which represents the removal (or outflow) of species $U$ from the system and  prevents unbounded growth of $U.$

    The second equation governs the spatio-temporal evolution of the concentration $v$ of species $V$ in a reacting and diffusing medium. The diffusion term $D_v\Delta v$  accounts for the spatial spreading of species $V$. The nonlinear coupling term $-u^2v$ models an autocatalytic reaction, leading to the consumption of species $V$ and coupling the dynamics of the two species. The term $B(1-v)$ represents a constant feed of species $V$ from an external reservoir and acts as a replenishment mechanism that drives the concentration toward the homogeneous equilibrium state $v=1$, thereby preventing depletion of $V.$  

    Let $\Omega \subset \mathbb{R}^d,$ with $d\in\{1,2,3\},$ be a bounded connected domain with sufficiently smooth boundary $\partial \Omega,$ and let $\nu(x)$ denote the outer unit normal vector at $x\in \partial\Omega.$ We use the notation $\partial_t\equiv \frac{\partial}{\partial t},~\partial_{\nu}\equiv \frac{\partial}{\partial \nu}.$ Throughout this paper, we assume that $D_u=D_v=1$ and that $A,B>0$ are constants in \eqref{Intro_Gray_Scott_System}. Consider the Gray-Scott system \eqref{Intro_Gray_Scott_System} with homogeneous Neumann boundary conditions:
    \begin{equation}
    \label{eq:uncontrolled-GS}
        \begin{cases}{}
			  \partial_{t}u - \Delta u -u^2v+Au = 0,\quad \text{in}\ (0,\infty)\times \Omega, \\
			  \partial_{t}v - \Delta v+u^2v-B(1-v) = 0, \quad \text{in}\ (0,\infty)\times \Omega, \\
			  {\partial_\nu} u = 0,\  {\partial_\nu} v = 0, \quad  \text{on}\ (0,\infty)\times \partial\Omega, \\
			  u(0,x)=u_0(x), \
			  v(0,x)=v_0(x), \quad x\in \Omega,
		\end{cases}
    \end{equation}

    To define the steady states of the system \eqref{eq:uncontrolled-GS}, we first introduce the following unbounded linear operator on $L^2(\Omega)$. The operator $(\Delta_N,\mathcal{D}(\Delta_N))$, defined by
    \begin{equation*}
        \mathcal{D}(\Delta_N) := \left\{ w\in H^2(\Omega): \partial_\nu w=0 \text{ on }\partial\Omega \right\}, \qquad \Delta_N w=\Delta w,
    \end{equation*}
    is called the Neumann Laplacian. Denote $\mathcal{D}(\Delta_N)^{2} \equiv \mathcal{D}(\Delta_N) \times \mathcal{D}(\Delta_N).$

    \begin{definition}(Steady-state).
		An element  $({\overline u},{\overline v})^{\top}\in \mathcal{D}(\Delta_N)^{2}$ is called a steady state of \eqref{eq:uncontrolled-GS} if
		  \begin{equation}
        \label{eq:stationary} 
			  \begin{cases}{}
				\Delta {\overline u} + \overline u^2 \overline v - A \overline u = 0,\quad x\in \Omega,\\
				  \Delta{ \overline v} - \overline u^2 \overline v + B(1- \overline v) = 0,\quad x\in \Omega,\\
				\partial_{\nu}\overline u = 0,~ \partial_{\nu} \overline v = 0,\quad x\in \partial\Omega.
			  \end{cases}         
		  \end{equation}
	\end{definition}
    Since $d\leq 3$, the Sobolev embedding
    \begin{equation*}
        H^2(\Omega)\hookrightarrow C(\overline{\Omega}) \hookrightarrow L^\infty(\Omega)
    \end{equation*}
    implies that $\overline u$ and $\overline v$ are continuous and bounded on $\overline{\Omega}$.
    
    For constant $(\overline u,\overline v)^{\top}$, system \eqref{eq:stationary} becomes  
    \begin{equation*}
        \overline{u}(\overline{u} \, \overline{v}-A)=0,~  -\overline{u}^2\overline{v}+B(1-\overline{v})=0.
    \end{equation*}
    Solving these equations, we obtain the trivial constant steady state $(0,1)^\top$ and, when $B\ge4A^2$, the nonzero constant steady states $(\overline u_+,\overline v_+)^\top$ and $(\overline u_-,\overline v_-)^\top$ given by
    \begin{align} 
    \label{eq:constant-equilibria}
        \begin{cases}
		      & (\overline{u}_+, \overline{v}_+)^\top =\bigg( \frac{\big(B+\sqrt{B(B-4A^2)}\big)}{2A},\ \frac{\big(B-\sqrt{B(B-4A^2)}\big)}{2B}\bigg)^\top, \\
		      &(\overline{u}_-, \overline{v}_-)^\top =\bigg(\frac{\big(B-\sqrt{B(B-4A^2)}\big)}{2A},\ \frac {\big(B+\sqrt{B(B-4A^2)}\big)}{2B}\bigg)^\top.
        \end{cases}     
    \end{align}
    When $B=4A^2$, the two nonzero steady states coincide, i.e., $({\overline{u}_+},{\overline{v}_+})^\top= ({\overline{u}_-}, {\overline{v}_-})^\top = ({B/2A},\ {1/2})^\top$. Moreover, for each nonzero constant steady state,
        $\overline u_+\overline v_+ = \overline u_-\overline v_- = A.$

    \begin{remark}
If a steady state $(\overline u,\overline v)^\top$ is
nonconstant, then both of its components must be nonconstant.
Indeed, if $\overline u$ is a nonzero constant, the first
equation of \eqref{eq:stationary} gives
$\overline v=A/\overline u$, which is constant.
If $\overline u=0$, the second equation, together with
the homogeneous Neumann boundary condition, yields
$\overline v=1$. Conversely, if $\overline v$ is constant,
the second equation of \eqref{eq:stationary} implies that
$\overline v\neq0$ and $\overline u^2$ is constant.
By continuity of $\overline u$ and connectedness of $\Om$,
it follows that $\overline u$ is constant.
\end{remark}

\subsection{The control problem}
    In this article, we are interested in the stabilization of steady states of the system \eqref{eq:uncontrolled-GS} by means of an interior control supported in a nonempty open subset of the domain $\Omega$ and acting only on the first equation. To this end, we consider the following controlled system corresponding to the system~\eqref{eq:uncontrolled-GS}
    \begin{equation}
    \label{eq:controlled-GS}
        \begin{cases}{}
			\partial_{t}u - \Delta u -u^2v+Au = \mathds{1}_{\omega} h,\quad \text{in}\ (0,\infty)\times \Omega, \\
			  \partial_{t}v - \Delta v+u^2v-B(1-v) = 0, \quad \text{in}\ (0,\infty)\times \Omega, \\
			  {\partial_\nu} u = 0,\  {\partial_\nu} v = 0, \quad  \text{on}\ (0,\infty)\times \partial\Omega, \\
			  u(0,x)=u_0(x), \
			  v(0,x)=v_0(x), \quad x\in \Omega,
		  \end{cases}
    \end{equation}
    where $\omega\subset\Omega$ is a nonempty open subset, $\mathds{1}_{\omega}$ denotes its characteristic function, and $h\in L^2((0,\infty)\times\Omega)$.  

    \begin{remark}
        Following the terminology used in \cite{BW03}, a steady state of the system \eqref{eq:uncontrolled-GS} will also be referred to as a steady state of the system \eqref{eq:controlled-GS}.
    \end{remark}
  

We now state the main result of this paper concerning
the local exponential stabilization of steady states
of the Gray-Scott system by finite-dimensional
internal feedback controls.

\begin{theorem}[Local exponential stabilization]
\label{thm:main}
Let $\omega\subset\Om$ be any nonempty open set and let
$(\overline u,\overline v)^\top\in D(\Delta_N)^2$ be a
steady state of \eqref{eq:uncontrolled-GS} with
$(\overline u,\overline v)^\top\neq(0,1)^\top$.
Then, for every prescribed decay rate $\gamma>0$,
there exist a finite-dimensional actuator space
\begin{equation*}
    \cU_\gamma
    =\operatorname{span}\{\psi_1,\ldots,\psi_m\}
    \subset L^2(\omega)
\end{equation*}
and a bounded feedback operator
$K_\gamma:L^2(\Om)\times L^2(\Om)\to\R^m$
such that the feedback \eqref{eq:feedback} locally
exponentially stabilizes $(\overline u,\overline v)^\top$
in $H^1(\Om)\times H^1(\Om)$.

More precisely, there exist $\rho,C>0$ such that,
for every
$(u_0,v_0)^\top\in H^1(\Om)\times H^1(\Om)$ satisfying
\begin{equation*}
    \norm{\binom{u_0-\overline u}{v_0-\overline v}}
    _{H^1(\Om)\times H^1(\Om)}<\rho,
\end{equation*}
the corresponding closed-loop system admits a unique
global mild solution satisfying
\begin{equation*}
    \binom{u-\overline u}{v-\overline v}
    \in C\big([0,\infty);H^1(\Om)\times H^1(\Om)\big)
\end{equation*}
and
\begin{equation*}
    \norm{\binom{u(t)-\overline u}{v(t)-\overline v}}
    _{H^1(\Om)\times H^1(\Om)}
    \le Ce^{-\gamma t}
    \norm{\binom{u_0-\overline u}{v_0-\overline v}}
    _{H^1(\Om)\times H^1(\Om)},
    \qquad t\ge0.
\end{equation*}
\end{theorem}

    The control and stabilization of semilinear parabolic systems is a challenging problem and has therefore attracted significant interest from researchers over the past two decades. Among the various control strategies, feedback stabilization, especially via finite-dimensional feedback controllers, is analytically more challenging but it has proven effective in stabilizing unstable steady states and is particularly relevant for practical applications. There has been extensive literature on results related to that; see, for instance, \cite{Bad12, Bar18, BW03, Liu10, YCB09}, and the references therein. Barbu and Wang \cite[Hypothesis~(H3)]{BW03} impose a structural condition requiring the off-diagonal derivatives of the reaction terms and a prescribed linear combination of the diagonal derivatives, evaluated at a steady state, to be spatially constant, with a nonzero coupling coefficient from the controlled to the uncontrolled component. Although the nontrivial constant steady states given in
\eqref{eq:constant-equilibria} satisfy this condition, linearization about a nonconstant steady state generally leads to a non-self-adjoint coupled operator with spatially varying coefficients, notably $2\overline u \, \overline v$ and $\overline u^2$. Consequently, a reduction to independent modes of the Neumann Laplacian is not available in general.

    Our approach accommodates this variable-coefficient, non-self-adjoint setting. The linearized operator is a bounded matrix-valued perturbation of the diagonal Neumann Laplacian and therefore generates an analytic semigroup with compact resolvent, leaving only finitely many critical spectral modes for any prescribed decay rate. These modes are isolated by a Riesz projection, and their controllability is established through the Fattorini-Hautus criterion, using unique continuation for the adjoint elliptic system and a nonvanishing property of the Gray-Scott steady-state coupling. This yields a finite-dimensional actuator space and a Riccati feedback that shifts the critical spectrum beyond the prescribed stability threshold. Analytic-semigroup smoothing and estimates for the quadratic-cubic nonlinear remainder then give local exponential stabilization of the full nonlinear system.

    Thus, the present work extends the finite-dimensional stabilization framework of \cite{BW03, YCB09} without imposing the structural hypothesis~(H3) of Barbu-Wang or the self-adjointness and a priori local linear-independence assumptions on unstable eigenfunctions employed by Yan, Coca and Barbu. In particular, the controllability required for finite-dimensional feedback is derived from the Gray-Scott PDE structure itself, rather than assumed as a spectral nondegeneracy condition. A useful structural observation is that, if a steady state $(\overline u,\overline v)\in \mathcal{D}(\Delta_N)^{2}$ satisfies
    \begin{equation*}
        \overline u \, \overline v\equiv0
    \end{equation*}
    on some nonempty open set, then the steady-states system~\eqref{eq:stationary} and scalar unique continuation imply
    \begin{equation*}
        (\overline u,\overline v)^{\top}=(0,1)^{\top}.
    \end{equation*}
    Consequently, for every other steady state and every nonempty open control region $\omega\subset\Omega$, the product $\overline u \, \overline v$ is nonzero at some point of $\omega$. By continuity, it is bounded away from zero on a smaller ball, which is precisely what is needed in the adjoint Hautus argument.

\subsection{Perturbation system and abstract formulation}
    Setting $y=u-\overline u, \, z=v-\overline v,$ and define the mapping 
    \begin{align}\label{eq:Nscalar}
        N:H^1(\Omega) \times H^1(\Omega) \longrightarrow L^2(\Omega), \,\, N(y,z) := \overline v\,y^2+2\overline u\,yz+y^2z.
    \end{align}
    Then the problem of stabilizing the steady state $(\overline u,\overline v)^\top$ of the system \eqref{eq:controlled-GS} is equivalent to the stabilizing zero steady state of the following system:
    \begin{equation}
    \label{eq:perturb}
        \begin{cases}
            y_t=\Delta y+(2\overline u \, \overline v-A)y+\overline u^2z+\mathds{1}_\omega h+N(y,z),\quad \text{in}\ (0,\infty)\times \Omega,\\
            z_t=\Delta z-2\overline u \, \overline v\,y-(\overline u^2+B)z-N(y,z),\quad \text{in}\ (0,\infty) \times \Omega,\\
            \partial_\nu y=\partial_\nu z=0, \quad \text{on}\ (0,\infty)\times\partial\Omega,\\
             y(0,x)=y_0(x):=u_0(x)-\overline u(x),\quad
             z(0,x)=z_0(x):=v_0(x)-\overline v(x), \quad x\in \Omega.
        \end{cases}
    \end{equation}
    We define the bounded linear control operator $\mathcal{B}:L^2(\omega)\to L^2(\Omega) \times L^2(\Omega)$ by
    \begin{equation}
    \label{eq:B}
        \mathcal{B} h=\begin{pmatrix}
                    \mathds{1}_\omega h\\
                    0
                    \end{pmatrix},
        \qquad h\in L^2(\omega),
    \end{equation}
    Let us denote $Y=(y,z)^\top,\, Y_0=(y_0,z_0)^\top.$ For $Y=(y,z)^\top$, define the unbounded linear operator $\cG:\mathcal{D}(\cG)\subset Z\to Z$ by
    \begin{equation}
    \label{eq:G}
        \mathcal{D}(\cG)=\mathcal{D}(\Delta_N)^2, \quad \mathcal{G} Y =
        \begin{pmatrix}
           \Delta+2\overline u\,\overline v-A & \overline u^2\\
           -2\overline u\,\overline v & \Delta-\overline u^2-B
        \end{pmatrix}Y.
    \end{equation}
    Define the nonlinear mapping
    \begin{equation*}
        \mathcal N:H^1(\Omega) \times H^1(\Omega) \longrightarrow L^2(\Omega) \times L^2(\Omega)
    \end{equation*}
    by
    \begin{equation}
    \label{eq:Nonlinear_map}
        \mathcal N(Y)=\binom{N(y,z)}{-N(y,z)},\qquad Y=(y,z)^\top\in H^1(\Omega) \times H^1(\Omega),
    \end{equation}
    Then system~\eqref{eq:perturb} can be written as the following evolution equation
    \begin{equation}
    \label{eq:evolution}
        \begin{cases}
            Y'(t)=\mathcal GY(t)+\mathcal Bh(t)+\mathcal N(Y(t)), & t>0,\\
            Y(0)=Y_0,
        \end{cases}
    \end{equation}
    in the Hilbert space $L^2(\Omega) \times L^2(\Omega).$

\subsection{Organization of the paper}
    The rest of the paper is organized as follows. Section~\ref{sec:functional} introduces the functional setting and notations. Section~\ref{sec:linear} establishes unique continuation, proves a nonvanishing property of the steady state system, and constructs a finite-dimensional feedback for the full variable-coefficient linearization. Section~\ref{sec:nonlinear} proves local exponential stabilization of the nonlinear Gray-Scott system. Section~\ref{sec:remarks} records consequences and special cases, and Section~\ref{sec:conclusion} concludes the paper.

\section{Functional setting and notation} 
\label{sec:functional}
    For $1\leq p\leq\infty$, we denote by $L^p(\Omega)$ the Lebesgue space equipped with the norm 
    \begin{align*}
        \|u\|_{L^p(\Omega)}=
        \begin{cases}
             \displaystyle\left(\int_\Omega |u(x)|^p\,dx\right)^{1/p},& 1\leq p<\infty,\\[3mm]
             \displaystyle\operatorname*{ess\,sup}_{x\in\Omega}|u(x)|,& p=\infty.
        \end{cases}
    \end{align*}

    For an integer $k\geq 0$ and $1\leq p\leq\infty$, the Sobolev space $W^{k,p}(\Omega)$ is defined by
    \begin{equation*}
        W^{k,p}(\Omega):= \left\{u \in L^p(\Omega) \; : \;D^{\alpha} u \in L^p(\Omega)\ \text{for all multi-indices } {\alpha} \text{ with } |{\alpha}|\le k\right\},
    \end{equation*}
    equipped with the norm
    \begin{align*}
        \|u\|_{W^{k,p}(\Omega)}:=
        \begin{cases}
            \Bigg( \displaystyle\sum_{|\alpha|\leq k}\|D^\alpha u\|^p_{L^p(\Omega)}\Bigg)^{1/p},& 1\leq p<\infty,\\[3mm]
            \displaystyle\max_{|\alpha|\leq k}\|D^\alpha u\|_{L^\infty(\Omega)},& p=\infty.
        \end{cases}
    \end{align*}
    We denote $H^k(\Omega): = W^{k,2}(\Omega),$ which is a Hilbert space, and the norm is induced by the standard inner product
    \begin{equation*}
        \langle u,v\rangle_{H^k(\Omega)}:= \sum_{|\alpha|\le k} \int_\Omega D^\alpha u(x)\, D^\alpha v(x)\, dx.
    \end{equation*}
    
    Let us denote
    \begin{equation*}
        H:=L^2(\Om),\qquad Z:=H\times H,\qquad X:=H^1(\Om)\times H^1(\Om),
    \end{equation*}
    where these spaces are real unless explicitly complexified. For the spectral arguments, we write
    \begin{equation*}
        Z_\C:=L^2(\Om;\C) \times L^2(\Om;\C) \equiv L^2(\Om;\C)^2
    \end{equation*}
    for the complexification of $Z,$ and use the same symbols $\cG$ and $\cB$ for the complex-linear extensions of the corresponding real operators.

\section{Finite-dimensional stabilization of the linearized system}
\label{sec:linear}
    Throughout Sections~\ref{sec:linear} and~\ref{sec:nonlinear}, we fix a steady state $(\overline u,\overline v)^\top$ of \eqref{eq:uncontrolled-GS} such that
    \begin{equation}
    \label{eq:Nontrivial_ss}
        (\overline u,\overline v)^\top\neq(0,1)^\top.
    \end{equation}
    Consider the linearized controlled system corresponding to \eqref{eq:evolution}:
    \begin{equation}
    \label{eq:linear}
        \begin{cases}
            Y_t=\mathcal GY+\mathcal Bh, & t>0,\\
            Y(0)=Y_0.
        \end{cases}
    \end{equation}

\begin{lemma}
\label{lem:compact}
    The operator $\mathcal G$ defined in \eqref{eq:G} generates an analytic semigroup on $Z$ and has compact resolvent. Moreover, the spectrum of its complexification on $Z_{\mathbb C}$ consists of isolated eigenvalues of finite algebraic multiplicity, with no finite accumulation point. For every $\eta>0$, the set
    \begin{equation*}
        \{\lambda\in\sigma(\mathcal G):\operatorname{Re}\lambda\ge -\eta\}
    \end{equation*}
    is finite.
\end{lemma}
\begin{proof}
    Write $\mathcal{G}=\mathcal{G}_0 + \mathcal{M},$ where
    \begin{equation*}
        \mathcal D(\mathcal G_0)=\mathcal \mathcal{D}(\Delta_N)^2,\qquad
        \mathcal G_0:=\begin{pmatrix}
                      \Delta_N & 0 \\
                       0 & \Delta_N
                       \end{pmatrix},\qquad
        \mathcal{M} = \begin{pmatrix}
                       2\overline{u} \, \overline{v}-A & \overline{u}^2\\
                      -2\overline{u} \, \overline{v} & -\overline{u}^2-B
                      \end{pmatrix}. 
    \end{equation*}
    Since $\overline u,\overline v\in C(\overline{\Omega})$, we have $\mathcal M\in L^\infty(\Omega;\mathbb R^{2\times2}).$ Hence, $\mathcal M$ induces a bounded multiplication operator on $Z$. Since $\mathcal G_0$ generates an analytic semigroup on $Z$, the bounded perturbation theorem implies that $\mathcal G=\mathcal G_0+\mathcal M$ also generates an analytic semigroup on $Z$.

    Moreover, the embedding
    \begin{equation*}
        \mathcal D(\mathcal G_0)=\mathcal \mathcal{D}(\Delta_N)^2 \hookrightarrow Z
    \end{equation*}
    is compact. Hence, $\mathcal G_0$ has compact resolvent. Since $\mathcal M$ is bounded, compactness of the resolvent is preserved under bounded perturbations, and therefore $\mathcal G$ also has compact resolvent. Let $\mathcal{G} Y=\lambda Y, \quad \norm{Y}_{Z}=1$. Then $\lambda = -\norm{\nabla Y}_{L^2}^2 + \ip{\cM Y}{Y}_Z$.
    Hence
    \begin{equation*}
        \operatorname{Re}\lambda \leq \norm{\mathcal{M}}_{\mathcal L(Z)},\qquad |\operatorname{Im} \lambda| \leq \norm{\mathcal{M}}_{\mathcal{L}(Z)}. 
    \end{equation*}
    Thus, the eigenvalues satisfying $\operatorname{Re}\lambda\geq-\eta$ lie in a bounded rectangle. Since $\mathcal G$ has compact resolvent, its spectrum consists of isolated eigenvalues of finite algebraic multiplicity, with no finite accumulation point (see, e.g., \cite[Chap.~III, \S6.8]{Kat95}). Consequently, only finitely many eigenvalues satisfy $\operatorname{Re}\lambda\geq-\eta$.  This completes the proof of Lemma~\ref{lem:compact}.
\end{proof}

\begin{lemma}[Unique continuation]
\label{lem:UCP}
    Let $\mathcal O\subset\R^d$ be connected, $1\le d\le3$, and let $P\in H^2_{\mathrm{loc}}(\mathcal O;\C^m)$ satisfy
    \begin{equation*}
        -\Delta P+V(x)P=0 \qquad\text{in }\mathcal O,
    \end{equation*}
    where $V\in L^\infty(\mathcal O;\C^{m\times m})$. If $P$ vanishes on a nonempty open subset of $\mathcal O$, then $P\equiv0$ in $\mathcal O$.
\end{lemma}
\begin{proof}
    Assume first that $d=1$. Then the connected set $\mathcal O$ is an interval and
    \begin{equation*}
        P''=VP \qquad\text{a.e. in }\mathcal O.
    \end{equation*}
    Since $P\in H^2_{\mathrm{loc}}(\mathcal O)$, both $P$ and $P'$ have continuous representatives on compact subintervals. If $P$ vanishes on a nonempty open interval, choose an interior point $x_0$ of that interval. Then $P(x_0)=P'(x_0)=0$. On every compact subinterval $I\Subset\mathcal O$ containing $x_0$, the vector $W=(P,P')^\top$ satisfies
    \begin{equation*}
        W'=\begin{pmatrix}
              0&I_m\\
              V&0
           \end{pmatrix}W
        \qquad\text{a.e. on }I,
    \end{equation*}
    with an $L^\infty$ coefficient matrix. Uniqueness for this first-order system (equivalently, Gronwall's inequality applied to its integral form) gives $W\equiv0$ on $I$. By connectedness, $P\equiv0$ in $\mathcal O$.

    Assume now that $d\ge2$. Write $P=P_1+iP_2$ and $V=V_1+iV_2$, and set
    \begin{equation*}
        U:=\binom{P_1}{P_2}\in H^2_{\mathrm{loc}}(\mathcal O;\R^{2m}).
    \end{equation*}
    Then
    \begin{equation*}
        -\Delta U+\widetilde V(x)U=0, \qquad
        \widetilde V=\begin{pmatrix}
                      V_1&-V_2\\
                      V_2&V_1
                      \end{pmatrix}
        \in L^\infty(\mathcal O;\R^{2m\times2m}).
    \end{equation*}
    Hence, for every component $U_k$,
\begin{equation*}
    |\Delta U_k|^2
    \le C\sum_{j=1}^{2m}|U_j|^2
    \le C\sum_{j=1}^{2m}
    \bigl(|\nabla U_j|^2+|U_j|^2\bigr)
    \qquad\text{a.e. in }\mathcal O.
\end{equation*}
Aronszajn's Remark~3 states that his main
unique-continuation theorem extends to a finite system
$(u_1,\ldots,u_N)$ satisfying inequalities of the form
\begin{equation*}
    |\mathcal L u_k|^2
    \le M\sum_{i=1}^N
    \left(
        \sum_{j=1}^d
        \left|\frac{\partial u_i}{\partial x^j}\right|^2
        +|u_i|^2
    \right),
    \qquad k=1,\ldots,N,
\end{equation*}
provided the same second-order elliptic operator
$\mathcal L$ acts on every component; see
\cite[p.~248, Remark~3]{Aro57}.
Taking $\mathcal L=\Delta$ and $N=2m$, the estimate
obtained above satisfies this condition.
Since $U$ vanishes on a nonempty open subset of
$\mathcal O$, it has a zero of infinite order at
any interior point of that subset.
Aronszajn's theorem therefore yields
$U\equiv0$ in the connected domain $\mathcal O$,
and hence $P\equiv0$.
\end{proof}
\begin{lemma} 
\label{lem:nonvanishing}
    Let $(\overline u,\overline v)^{\top} \in \mathcal{D}(\Delta_N)^2$ be a steady state satisfying \eqref{eq:Nontrivial_ss}. Then, for every nonempty open set $\omega\subset\Om$,
    \begin{equation*}
        \overline u \, \overline v\not\equiv0 \qquad\text{in }\omega.
    \end{equation*}
    Equivalently, every nonempty open $\omega\subset\Om$ contains a point $x_0$ such that $\overline u(x_0)\overline v(x_0)\neq0$.
\end{lemma}
\begin{proof}
    Suppose, to the contrary, that
    \begin{equation*}
        \overline u \, \overline v=0 \quad\text{in some nonempty open set }\omega\subset\Om.
    \end{equation*}
    We first show that $\overline u\equiv0$ in $\omega$. If $\overline u(x_0)\neq0$ at some $x_0\in\omega$, continuity of $\overline u$ gives $r>0$ such that
    \begin{equation*}
        B_r(x_0)\Subset\omega, \qquad \overline u(x)\neq0\quad (\forall x\in B_r(x_0)).
    \end{equation*}
    Since $\overline u \, \overline v=0$ in $\omega$, it follows that $\overline v=0$ in $B_r(x_0)$. Hence $\Delta\overline v=0$ almost everywhere in $B_r(x_0)$. The second equation of \eqref{eq:stationary} then yields $B=0,$ contradicting $B>0$. Thus
    \begin{equation*}
        \overline u\equiv0 \quad \text{in}\,\, \omega. 
    \end{equation*}
    The first equation of \eqref{eq:stationary} can be rewritten as
    \begin{equation*}
        -\Delta\overline u+(A-\overline u \, \overline v)\overline u=0 \qquad\text{in }\Om.
    \end{equation*}
    Because $\overline u \, \overline v\in L^\infty(\Om)$, Lemma~\ref{lem:UCP} with $m=1$ implies from $\overline u=0$ on the nonempty open set $\omega$ that
    \begin{equation*}
        \overline u\equiv0 \qquad\text{in }\Om.
    \end{equation*}
    The second equation of \eqref{eq:stationary} therefore becomes
    \begin{equation*}
        -\Delta\overline v+B\overline v=B, \qquad \partial_\nu\overline v=0.
    \end{equation*}
    Setting $w=\overline v-1$, we obtain
    \begin{equation*}
        -\Delta w+Bw=0,\qquad \partial_\nu w=0.
    \end{equation*}
    Multiplying the equation by $w$ and integrating over $\Om$, we obtain
    \begin{equation*}
        \int_\Om |\nabla w|^2\,dx+B\int_\Om |w|^2\,dx=0,
    \end{equation*}
    so $w\equiv0$. Hence $(\overline u,\overline v)^{\top}=(0,1)^{\top}$, a contradiction.
\end{proof}

\begin{remark}
    For every steady state satisfying \eqref{eq:Nontrivial_ss}, the set $\{x\in\Om:\overline u(x)\overline v(x)\neq0\}$ is open and dense in $\Om$. Indeed, openness follows from continuity, while density follows from Lemma~\ref{lem:nonvanishing}.
\end{remark}

    We next verify the Fattorini-Hautus condition. To this end, we introduce the adjoint operators $\mathcal G^*$ and $\mathcal B^*$. The adjoint operator $\mathcal{G}^*:\mathcal D(\mathcal{G}^*)\subset Z\to Z$ of $\mathcal{G}$ is given by
    \begin{equation}
    \label{eq:Gstar}
        \mathcal D(\mathcal{G}^*)=\mathcal \mathcal{D}(\Delta_N)^2, \quad \mathcal{G}^*=\begin{pmatrix}
                       \Delta_N+2\overline u\,\overline v-A & -2\overline u\,\overline v\\
                       \overline u^2 & \Delta_N-\overline u^2-B
                      \end{pmatrix}.
    \end{equation}

    The adjoint operator $\cB^*:Z\to L^2(\omega)$ of $\mathcal{B}$ is given by
    \begin{equation}
    \label{eq:Bstar}
        \cB^*\begin{pmatrix} p\\q \end{pmatrix}
        = \mathds{1}_{\omega} p =p|_\omega, \qquad (p,q)^\top\in Z.
    \end{equation}

\begin{proposition}[Fattorini-Hautus condition]
\label{prop:Hautus}
    For every $\lambda\in\C$,
    \begin{equation*}
        \ker(\lambda I-\cG^*) \cap \ker\cB^* = \{0\}.
    \end{equation*}
\end{proposition}
\begin{proof}
    Let $P=(p,q)^\top\in D(\cG^*)$ satisfy 
    $\mathcal{G}^*P=\lambda P$ and $\cB^*P=0$. By \eqref{eq:Bstar}, $p=0$ in $\omega$. Since $p\in H^2(\Om)$, it follows that $\Delta p=0$ in $\omega.$ The first equation of $\cG^*P=\lambda P$ is
    \begin{equation*}
        \Delta p+(2\overline u \, \overline v-A)p-2\overline u \, \overline v\,q=\lambda p.
    \end{equation*}
    Therefore $2\overline u \, \overline v\,q=0$ in $\omega$. By Lemma~\ref{lem:nonvanishing}, there exists $x_0\in\omega$ such that $\overline u(x_0)\overline v(x_0)\ne0$. By continuity, there exist $r>0$ and $c_0>0$ such that $B_r(x_0)\Subset\omega$ and
    \begin{equation*}
        |\overline u(x)\overline v(x)|\ge c_0 \qquad (x\in B_r(x_0)).
    \end{equation*}
     Hence $q=0$ in $B_r(x_0)$, and therefore $p=q=0$ in $B_r(x_0)$.

    The eigenvalue equation $ \mathcal{G}^*P=\lambda P $ is an elliptic system $-\Delta P+V_\lambda(x)P=0$ with $V_\lambda\in L^\infty(\Om;\C^{2\times2})$. Lemma~\ref{lem:UCP} gives $P\equiv0$ in $\Om$.
\end{proof}

    Fix $\gamma>0$ and define
    \begin{equation}
    \label{eq:critical}
        \Sigma_\gamma:= \{\lambda\in\sigma(\cG): \operatorname{Re} \lambda\ge-4\gamma\}.
    \end{equation}
    By Lemma~\ref{lem:compact}, $\Sigma_\gamma$ is finite. Since $\cG$ has real coefficients, $\lambda\in\sigma(\cG)\Longrightarrow \overline\lambda\in\sigma(\cG),$
    and therefore
    \begin{equation}
    \label{eq:sigma-symmetry}
        \overline{\Sigma_\gamma}=\Sigma_\gamma.
    \end{equation}
    If $\Sigma_\gamma$ is empty, the linearized system \eqref{eq:linear} with $h=0$ is already exponentially stable with the prescribed decay rate $\gamma$, and no feedback is needed. Assume $\Sigma_\gamma\ne\varnothing$.

    Choose $\delta>0$ such that
    \begin{equation*}
        \sigma(\cG)\cap \{\lambda\in\C:-4\gamma-\delta\le\operatorname{Re} \lambda<-4\gamma\} =\varnothing.
    \end{equation*}
    Such a $\delta$ exists because $\Sigma_\gamma$ is finite and the spectrum has no finite accumulation point. Using the spectral bounds in Lemma~\ref{lem:compact}, one may take $\Gamma$ to be a positively oriented rectangle with left side $\operatorname{Re} \lambda=-4\gamma-\delta/2$ and with the remaining three sides chosen outside the spectral rectangle. We choose it symmetric with respect to the real axis, so that $\Gamma=\overline{\Gamma}$ as a set and $\operatorname{ind}_\Gamma(\overline\lambda)= \operatorname{ind}_\Gamma(\lambda)$. Then $\Gamma$ encloses precisely $\Sigma_\gamma$. Define
    \begin{equation}
    \label{eq:Riesz}
        \Pi_\gamma := \frac{1}{2\pi i} \int_\Gamma (\lambda I-\cG)^{-1}\,d\lambda.
    \end{equation}
    Set
    \begin{equation*}
        Z_\gamma^u:= \operatorname{Ran}\Pi_\gamma,\qquad Z_\gamma^s:= \ker\Pi_\gamma. 
    \end{equation*}
    By the standard Riesz-projection theory for isolated spectral sets,
    \begin{equation*}
        Z_\C=Z_\gamma^u\oplus Z_\gamma^s,\qquad \dim_{\C}Z_\gamma^u<\infty,  
    \end{equation*}
    and both subspaces are $\cG$-invariant; see, for example, \cite[Chap.~III, \S6.4]{Kat95}. Because $\Gamma$ is conjugation-symmetric and $\cG$ is real, $\Pi_\gamma$ commutes with complex conjugation.

    Define
    \begin{align}
    \label{eq:actuator-space}
        \cU_{\gamma,\C}= \operatorname{span}_{\C} \left\{ \cB^*\Phi:\Phi\in\ker(\lambda I-\cG^*),\ \lambda\in\Sigma_\gamma\right\} \subset L^2(\omega;\C).
    \end{align}
    The space is finite dimensional. It is also invariant under complex conjugation: if $\cG^*\Phi=\lambda\Phi$ with $\lambda\in\Sigma_\gamma$, then $\cG^*\overline\Phi=\overline\lambda\,\overline\Phi$, and \eqref{eq:sigma-symmetry} together with the reality of $\cB$ gives $\overline{\cB^*\Phi}=\cB^*\overline\Phi\in\cU_{\gamma,\C}$. Hence there exists a finite-dimensional real space $\cU_\gamma\subset L^2(\omega;\R)$ whose complexification is $\cU_{\gamma,\C}$. Let $m=\dim\cU_\gamma$ and choose a real orthonormal basis $\psi_1,\ldots,\psi_m$. Define the finite-dimensional actuator injection $\mathcal E_\gamma:\R^m\to Z$ by
    \begin{equation}
    \label{eq:actuator-injection}
        \mathcal E_\gamma a = \bigg( \mathds{1}_{\omega} \sum_{j=1}^m a_j\psi_j,\, 0 \bigg)^{\top}.
    \end{equation}
    We use the same symbol for its complex-linear extension $\mathcal E_\gamma:\C^m\to Z_\C$ when applying the complex PBH criterion. Then
    \begin{equation*}
        \mathcal E_\gamma^*\binom pq=\left(\ip{p}{\psi_1}_{L^2(\omega)},\ldots, \ip{p}{\psi_m}_{L^2(\omega)}\right)^\top.
    \end{equation*}
    Let $Z_{\gamma,\R}^u=\{Y\in Z_\gamma^u:\overline Y=Y\}$ be the real critical subspace. Since $\Pi_\gamma$ commutes with complex conjugation,
    \begin{equation*}
        Z_\gamma^u=Z_{\gamma,\R}^u\oplus iZ_{\gamma,\R}^u, \qquad \dim_{\C}Z_\gamma^u = \dim_{\R}Z_{\gamma,\R}^u.
    \end{equation*}
    Set $n_\gamma:=\dim_{\R}Z_{\gamma,\R}^u$. Choose a real basis of $Z_{\gamma,\R}^u$ and let
    \begin{equation*}
        T_\gamma:\R^{n_\gamma}\longrightarrow Z_{\gamma,\R}^u
    \end{equation*}
    be the corresponding coordinate isomorphism. We use the same symbol for its complex-linear extension
    \begin{equation*}
        T_\gamma:\C^{n_\gamma}\longrightarrow Z_\gamma^u.
    \end{equation*}
    Define the real matrices $\mathsf G_\gamma:=T_\gamma^{-1}\cG T_\gamma$ and $\mathsf B_\gamma:=T_\gamma^{-1}\Pi_\gamma\mathcal E_\gamma$.

\begin{proposition}
\label{prop:controllable}
    The finite-dimensional real pair $(\mathsf G_\gamma,\mathsf B_\gamma)$ is controllable.
\end{proposition}
\begin{proof}
    Controllability is invariant under a change of coordinates, and for real matrices the Kalman rank over $\R$ is the same as the Kalman rank of the complexified pair over $\C$. We therefore verify the complex PBH criterion. We use the convention that the complex Hilbert space inner product is linear in its first argument.

    Let $w\in\C^{n_\gamma}$ be a left eigenvector of $\mathsf G_\gamma$ such that $w^*\mathsf G_\gamma=\lambda w^*$ and $w^*\mathsf B_\gamma=0$. Define a continuous linear functional on $Z_\C$ by
    \begin{equation*}
        \ell(y):=w^*T_\gamma^{-1}\Pi_\gamma y.
    \end{equation*}
    By the Riesz representation theorem, there exists a unique $\Phi\in Z_\C$ such that
    \begin{equation*}
        \ell(y)=\langle y,\Phi\rangle_Z \qquad\text{for all }y\in Z_\C.
    \end{equation*}
    Since $\ell((I-\Pi_\gamma)y)=0$, one has $\Phi\in\operatorname{Ran}\Pi_\gamma^*$. Moreover, for every $y\in D(\cG)$, using $\Pi_\gamma\cG=\cG\Pi_\gamma$ and $\mathsf G_\gamma=T_\gamma^{-1}\cG T_\gamma$, we obtain
    \begin{equation*}
        \langle \cG y,\Phi\rangle_Z =w^*T_\gamma^{-1}\Pi_\gamma\cG y=w^*\mathsf G_\gamma T_\gamma^{-1}\Pi_\gamma y =\lambda\langle y,\Phi\rangle_Z =\langle y,\overline\lambda\Phi\rangle_Z.
    \end{equation*}
    Thus the map $y\mapsto\langle\cG y,\Phi\rangle_Z$ is represented on $D(\cG)$ by the vector $\overline\lambda\Phi$. By the definition of the Hilbert space adjoint, this proves simultaneously that $\Phi\in D(\cG^*)$ and $\cG^*\Phi=\overline\lambda\Phi$. Here the conjugate eigenvalue is essential because the inner product is linear in its first argument. Likewise, for every $a\in\C^m$,
    \begin{equation*}
        \langle \mathcal E_\gamma a,\Phi\rangle_Z =w^*T_\gamma^{-1}\Pi_\gamma \mathcal E_\gamma a =w^*\mathsf B_\gamma a=0,    
    \end{equation*}
    so $\mathcal E_\gamma^*\Phi=0.$ Writing $\Phi=(p,q)^\top$, this means
    \begin{equation*}
        \langle p,\psi_j\rangle_{L^2(\omega)}=0,\qquad j=1,\ldots,m.
    \end{equation*}
    Since $\lambda$ is an eigenvalue of the critical matrix $\mathsf G_\gamma$, we have $\lambda\in\Sigma_\gamma$. Hence $\overline\lambda\in\Sigma_\gamma$ by \eqref{eq:sigma-symmetry}. Together with $\cG^*\Phi=\overline\lambda\Phi$, the definition \eqref{eq:actuator-space} shows that $\cB^*\Phi=p|_\omega$ belongs to $\cU_{\gamma,\C}$. It is also orthogonal to $\cU_{\gamma,\C}$ by $\mathcal E_\gamma^*\Phi=0$. Therefore $\cB^*\Phi=0$. Proposition~\ref{prop:Hautus} gives $\Phi=0$. Finally, taking $y=T_\gamma z$ with arbitrary $z\in\C^{n_\gamma}$ gives $w^*z=0$, hence $w=0$. Thus no nonzero left eigenvector can annihilate $\mathsf B_\gamma$, and the PBH criterion proves that $(\mathsf G_\gamma,\mathsf B_\gamma)$ is controllable.
\end{proof}

    Consider the matrix algebraic Riccati equation
    \begin{equation}
    \label{eq:ARE}
        (\mathsf G_\gamma+4\gamma I)^\top P_\gamma+ P_\gamma(\mathsf G_\gamma+4\gamma I) - P_\gamma\mathsf B_\gamma\mathsf B_\gamma^\top P_\gamma + I =0.
    \end{equation}
    Since $(\mathsf G_\gamma,\mathsf B_\gamma)$ is controllable and the state weight is the identity, the standard continuous-time algebraic Riccati theorem yields a unique positive definite stabilizing solution $P_\gamma$ of \eqref{eq:ARE}; see, for example, \cite{AM90}. Define $\mathsf K_\gamma^u:=-\mathsf B_\gamma^\top P_\gamma$ and the full-state feedback coefficient map
    \begin{equation}
    \label{eq:Kfull}
        K_\gamma:= \mathsf K_\gamma^uT_\gamma^{-1}\Pi_\gamma:Z\longrightarrow\R^m.
    \end{equation}
    The use of a real basis makes the feedback manifestly real-valued. In the preceding PBH argument, $\cG^*$ and $\Pi_\gamma^*$ are Hilbert space adjoints on the full complexified space; in \eqref{eq:ARE} the symbol $\top$ denotes the ordinary transpose of the finite-dimensional real matrices. This avoids any identification of $\operatorname{Ran}\Pi_\gamma^*$ with $\operatorname{Ran}\Pi_\gamma$.

    The distributed feedback is
    \begin{equation}
    \label{eq:feedback}
        h(t,x)=\sum_{j=1}^m [K_\gamma Y(t)]_j\psi_j(x).
    \end{equation}

\begin{theorem}[Linear stabilization]
\label{thm:linear-stab}
    The closed-loop operator
    \begin{equation*}
        \cG_{\rm cl}= \cG+\mathcal E_\gamma K_\gamma,\qquad D(\cG_{\rm cl})=D(\cG),
    \end{equation*}
    generates an analytic exponentially stable semigroup. In particular,
    \begin{equation*}
        \norm{e^{t\cG_{\rm cl}}}_{\mathcal L(Z)} \le M_\gamma e^{-3\gamma t}, \qquad t\ge0,
    \end{equation*}
    and, after reducing the rate if necessary,
    \begin{equation}
    \label{eq:smoothing}
        \norm{e^{t\cG_{\rm cl}}}_{\mathcal L(Z,X)} \le C_\gamma t^{-1/2}e^{-2\gamma t}, \qquad t>0,
    \end{equation}
    together with
    \begin{equation}
    \label{eq:X-semigroup}
        \norm{e^{t\cG_{\rm cl}}}_{\mathcal L(X)}\le C_\gamma e^{-2\gamma t}, \qquad t\ge0.
    \end{equation}
\end{theorem}
\begin{proof}
    Let $\cG_\gamma^s:=\cG|_{Z_\gamma^s}$. With respect to $Z_\C=Z_\gamma^u\oplus Z_\gamma^s$,
    \begin{equation*}
        \cG_{\rm cl}= \begin{pmatrix}
                      T_\gamma(\mathsf G_\gamma+\mathsf B_\gamma\mathsf K_\gamma^u)T_\gamma^{-1} & 0\\
                      (I-\Pi_\gamma)\mathcal E_\gamma\mathsf K_\gamma^uT_\gamma^{-1} & \cG_\gamma^s
                      \end{pmatrix}.
    \end{equation*}
    The Riccati construction gives $s(T_\gamma(\mathsf G_\gamma+\mathsf B_\gamma\mathsf K_\gamma^u)T_\gamma^{-1})<-4\gamma$. By definition of $\Sigma_\gamma$, we have $\sigma(\cG_\gamma^s)\subset\{\operatorname{Re} \lambda<-4\gamma\}$. Compactness of the resolvent provides a spectral gap from the line $\operatorname{Re} \lambda=-4\gamma$. Since analytic semigroups are spectrum determined, each diagonal semigroup decays with some rate strictly larger than $4\gamma$. The lower-left block is a convolution of the two diagonal semigroups, hence
    \begin{equation*}
        \norm{e^{t\cG_{\rm cl}}}\le C(1+t)e^{-4\gamma t} \le M_\gamma e^{-3\gamma t}.
    \end{equation*}
    Finally, $\cG_{\rm cl}$ differs from $\Delta_NI_2$ by a bounded operator. Therefore, after a sufficiently large shift $\eta>0$,
    \begin{equation*}
        D((\eta I-\cG_{\rm cl})^{1/2})=H^1(\Om)^2=X
    \end{equation*}
    with equivalent norms. Standard fractional-power estimates for analytic semigroups yield \eqref{eq:smoothing} and \eqref{eq:X-semigroup}; see \cite{Hen81, Paz83}.
\end{proof}

\section{Nonlinear stabilization}
\label{sec:nonlinear}

\begin{lemma}
\label{lem:N}
    The map $\mathcal{N}:X\to Z$ given by \eqref{eq:Nonlinear_map} is locally Lipschitz and satisfies
    \begin{equation*}
        \norm{\cN(Y)}_Z \le C(\norm{Y}_X^2+\norm{Y}_X^3).
    \end{equation*}
    Moreover, for every $R>0$,
    \begin{equation*}
        \norm{\cN(Y_1)-\cN(Y_2)}_Z \le C(R+R^2)\norm{Y_1-Y_2}_X
    \end{equation*}
    whenever $\norm{Y_i}_X\le R$.
\end{lemma}
\begin{proof}
    Since $d\le3$, we have $H^1(\Om)\hookrightarrow L^4(\Om)\cap L^6(\Om)$, and hence, by the definition of $N$ in \eqref{eq:Nscalar},
    \[\begin{aligned}
        \norm{N(y,z)}_2&\le \norm{\overline v}_\infty\norm{y}_4^2 + 2\norm{\overline u}_\infty\norm{y}_4\norm{z}_4 + \norm{y}_6^2\norm{z}_6\\
        &\le C(\norm{(y,z)}_X^2+\norm{(y,z)}_X^3).
    \end{aligned}\]
    For the difference estimate, use $y_1^2-y_2^2=(y_1-y_2)(y_1+y_2)$, $y_1z_1-y_2z_2=(y_1-y_2)z_1+y_2(z_1-z_2)$, and $y_1^2z_1-y_2^2z_2=(y_1^2-y_2^2)z_1+y_2^2(z_1-z_2)$, followed by Hölder's inequality and the same Sobolev embeddings.
\end{proof}

\begin{proof}[proof of Theorem~\ref{thm:main}]
    By \eqref{eq:feedback},
    \begin{equation*}
        Y_t=\cG_{\rm cl}Y+\cN(Y),\qquad Y(0)=Y_0.
    \end{equation*}
    Consider
    \begin{equation*}
        \cX_\gamma=\left\{Y\in C([0,\infty);X): \norm{Y}_{\cX_\gamma}= \sup_{t\ge0}e^{\gamma t}\norm{Y(t)}_X<\infty \right\}
    \end{equation*}
    and the ball $\mathbb B_R= \{Y\in\cX_\gamma:\norm{Y}_{\cX_\gamma}\le R\}$. Define
    \begin{equation*}
        (\mathcal TY)(t)=e^{t\cG_{\rm cl}}Y_0+\int_0^t e^{(t-s)\cG_{\rm cl}}\cN(Y(s))\,ds.
    \end{equation*}
    For $Y\in\mathbb B_R$, the semigroup term is continuous in $X$. Moreover, by the smoothing estimate and the local continuity of $\cN: X \to Z$, the convolution term is also continuous in $X$; hence $\mathcal T Y \in C([0,\infty); X)$. If $Y\in\mathbb B_R$, then $\norm{Y(s)}_X\le Re^{-\gamma s}$ and Lemma~\ref{lem:N} gives
    \begin{equation*}
        \norm{\cN(Y(s))}_Z \le C(R^2+R^3)e^{-\gamma s}.
    \end{equation*}
    By Theorem~\ref{thm:linear-stab},
    \begin{equation*}
        e^{\gamma t}\norm{(\mathcal TY)(t)}_X \le C_\gamma\norm{Y_0}_X+ C_\gamma(R^2+R^3) \int_0^t (t-s)^{-1/2}e^{-\gamma(t-s)}\,ds.
    \end{equation*}
    Since
    \begin{equation*}
        \int_0^\infty r^{-1/2}e^{-\gamma r}\,dr<\infty,
    \end{equation*}
    there are $C_0,C_1>0$ such that
    \begin{equation*}
        \norm{\mathcal TY}_{\cX_\gamma} \le C_0\norm{Y_0}_X+C_1(R^2+R^3).
    \end{equation*}
    Similarly, $\norm{\mathcal TY_1-\mathcal TY_2}_{\cX_\gamma} \le C_2(R+R^2) \norm{Y_1-Y_2}_{\cX_\gamma}$. Choose $R>0$ so small that
    \begin{equation*}
        C_2(R+R^2)\le\frac12 \qquad\text{and}\qquad C_1(R+R^2)\le\frac12.
    \end{equation*}
    Then $C_1(R^2+R^3)\le R/2$. Choose $\rho>0$ so that $C_0\rho\le\frac{R}{2}$. It follows that $C_0\rho+C_1(R^2+R^3)\le R$. The Banach fixed-point theorem gives a unique fixed point $Y\in\mathbb B_R$, hence a global mild solution. Since $\cN(0)=0$, the same Lipschitz estimate used above gives, for this fixed point,
    \begin{equation*}
        \norm{Y}_{\cX_\gamma}\le C_0\norm{Y_0}_X +C_2(R+R^2)\norm{Y}_{\cX_\gamma} \le C_0\norm{Y_0}_X+\frac12\norm{Y}_{\cX_\gamma}.
    \end{equation*}
    Therefore $\norm{Y}_{\cX_\gamma} \le 2C_0\norm{Y_0}_X$, and hence
    \begin{equation*}
        \norm{Y(t)}_X \le 2C_0e^{-\gamma t}\norm{Y_0}_X, \qquad t\ge0.
    \end{equation*}
    To obtain uniqueness among all mild solutions in $C([0, T]; X)$, let $Y_1, Y_2$ have the same initial datum. On any finite interval they lie in a common bounded subset of $X$, so the smoothing estimate and Lemma~\ref{lem:N} give
    \begin{equation*}
        \norm{Y_1(t)-Y_2(t)}_X \le C_T\int_0^t (t-s)^{-1/2} \norm{Y_1(s)-Y_2(s)}_X\,ds.
    \end{equation*}
    Henry's singular Gronwall inequality \cite[Lemma~7.1.1]{Hen81} yields $Y_1=Y_2$ on $[0,T]$.
\end{proof}

\section{Consequences and special cases}
\label{sec:remarks}

\begin{remark}
    Theorem~\ref{thm:main} holds for every nonempty open set $\omega\subset\Om$. The control region is arbitrary but fixed when the feedback system is defined. The actuator functions $\psi_j$ and the feedback map $K_\gamma$ are constructed after $\omega$, $\gamma$, and the steady state $(\overline u,\overline v)^{\top}$ have been fixed, and in general depend on all three.
\end{remark}

\begin{remark}  
The conclusion of Theorem~\ref{thm:main} remains valid
if the control acts only on the second equation.
In this case, the control operator
$\mathcal B_2:L^2(\omega)\to Z$ and its adjoint
are given by
\begin{equation*}
    \mathcal B_2h=\binom{0}{\mathds{1}_\omega h},
    \qquad
    \mathcal B_2^*\binom{p}{q}=q|_\omega.
\end{equation*}
 If $\mathcal{G}^*P=\lambda P$ and $\mathcal{B}_{2}^* P=0$, where $P=(p,q)^\top$, then $q=0$ in $\omega$. The second equation of $\mathcal{G}^*P=\lambda P$ therefore yields $\overline u^2p=0$ in $\omega$. By Lemma~\ref{lem:nonvanishing}, $\overline u^2\not\equiv0$ in every nonempty open subset of $\Om$. Hence, arguing as in the proof of Proposition~\ref{prop:Hautus}, we obtain the Fattorini-Hautus condition. The rest of the proof remains unchanged.
\end{remark}

 \begin{proposition}[Stabilization of the steady state $(0,1)^{\top}$]
    The following assertions hold for the trivial steady state $(0,1)^\top$.
    \begin{enumerate}
        \item It is locally exponentially stable in $X$ without control, with maximal decay rate $\min\{A,B\}$.
        \item It is locally exponentially stabilizable in $X$ by a finite-dimensional feedback control acting only on the first equation, with maximal decay rate $B$.
        \item It is locally exponentially stabilizable in $X$ by a finite-dimensional feedback control acting only on the second equation, with maximal decay rate $A$.
        \item It is locally exponentially stabilizable in $X$ at any prescribed decay rate by finite-dimensional feedback controls acting on both equations.
    \end{enumerate}
\end{proposition}

\begin{proof}
    Setting $y=u$ and $z=v-1$, we obtain the following perturbation system:
    \begin{equation*}
    \begin{cases}
        y_t=\Delta y-Ay+y^2+y^2z+\mathds{1}_\omega h_1,
        \quad \text{in } \, (0,\infty)\times\Om,\\
        z_t=\Delta z-Bz-y^2-y^2z+\mathds{1}_\omega h_2,
        \quad \text{in } \, (0,\infty)\times\Om,\\
        \partial_\nu y=\partial_\nu z=0,
        \quad \text{on } \,(0,\infty)\times\partial\Om,\\
        y(0,x)=y_0(x):=u_0(x),\ 
        z(0,x)=z_0(x):=v_0(x)-1, \quad
        x\in\Om,
    \end{cases}
    \end{equation*}
    where $h_1$ and $h_2$ denote the controls acting on the first and second equations, respectively. The corresponding linearized controlled system is
    \begin{equation*}
    \begin{cases}
        y_t=\Delta y-Ay+\mathds{1}_\omega h_1,
        &\text{in }(0,\infty)\times\Om,\\
        z_t=\Delta z-Bz+\mathds{1}_\omega h_2,
        &\text{in }(0,\infty)\times\Om,\\
        \partial_\nu y=\partial_\nu z=0,
        &\text{on }(0,\infty)\times\partial\Om,\\
        y(0,x)=y_0(x),\quad z(0,x)=z_0(x),
        &x\in\Om.
    \end{cases}
    \end{equation*}
    We first study the exponential stability of the linearized system without control. Setting $h_1=h_2=0$, the corresponding linear operator is given by
    \begin{equation*}
        \widetilde{\mathcal A}_0=\begin{pmatrix}
                                \Delta_N-A & 0\\
                                 0 & \Delta_N-B
                                  \end{pmatrix},
        \qquad D(\widetilde{\mathcal A}_0)=\mathcal{D}(\Delta_N)^2.
    \end{equation*}
    Since $\Delta_N$ generates an analytic contraction semigroup on $L^2(\Om)$, the operator $\widetilde{\mathcal A}_0$ generates an analytic semigroup on $Z$ given by
    \begin{equation*}
        e^{t\widetilde{\mathcal A}_0}=\begin{pmatrix}
        e^{-At}e^{t\Delta_N} & 0\\
        0 & e^{-Bt}e^{t\Delta_N}
    \end{pmatrix}.
    \end{equation*}
    Moreover, since $e^{t\Delta_N}$ is contractive on both $L^2(\Om)$ and $H^1(\Om)$, we obtain
    \begin{equation*}
    \max\left\{
        \norm{e^{t\widetilde{\mathcal A}_0}}_{\mathcal L(Z)},
        \norm{e^{t\widetilde{\mathcal A}_0}}_{\mathcal L(X)}
    \right\}
    \le e^{-\min\{A,B\}t}, \qquad t\ge0.
    \end{equation*}
    For $c>0$, consider the scalar controlled equation
    \begin{equation*}
    \begin{cases}
        w_t=\Delta w-cw+\mathds{1}_\omega h,
        &\text{in }(0,\infty)\times\Om,\\
        \partial_\nu w=0,
        &\text{on }(0,\infty)\times\partial\Om.
    \end{cases}
    \end{equation*}
    Since $\Delta_N-cI$ is self-adjoint and the adjoint of the control operator $h\mapsto\mathds{1}_\omega h$ is given by $p\mapsto p|_\omega$, the corresponding Fattorini-Hautus condition follows from Lemma~\ref{lem:UCP}. Indeed, if $p\in \mathcal{D}(\Delta_N)$ in the complexified space satisfies
    \begin{equation*}
        (\Delta_N-cI)p=\lambda p, \qquad p|_\omega=0,
    \end{equation*}
    for some $\lambda\in\mathbb{C}$, then Lemma~\ref{lem:UCP} yields $p\equiv0$ in $\Om$. Consequently, the finite-dimensional feedback construction of Section~3 applies to this scalar equation, yielding exponential stabilization at any prescribed decay rate.
    
    Thus, when the control acts only on the first equation, we can stabilize the first component of the linearized system with a decay rate strictly larger than $B$, while the second component decays with rate $B$. Similarly, when the control acts only on the second equation, we can stabilize the second component with a decay rate strictly larger than $A$, while the first component decays with rate $A$. When controls act on both equations, the above scalar feedback construction can be applied independently to each equation of the linearized system, yielding exponential stabilization at any prescribed decay rate.

    In each of the above cases, let $\mathcal A_{\mathrm{cl}}$ denote the corresponding linear closed-loop generator, with $\mathcal A_{\mathrm{cl}}=\widetilde{\mathcal A}_0$ in the uncontrolled case. Let $\gamma$ denote $\min\{A, B\}$, $B$, and $A$ in the first three cases, respectively, and any prescribed positive decay rate in the fourth case. By the scalar semigroup estimates and the smoothing estimates analogous to those in Theorem~\ref{thm:linear-stab}, for $t>0$, we obtain
    \begin{equation*}
         \norm{e^{t\mathcal A_{\mathrm{cl}}}}_{\mathcal L(X)}
         \le Ce^{-\gamma t}\quad \text{ and } \quad 
         \norm{e^{t\mathcal A_{\mathrm{cl}}}}_{\mathcal L(Z,X)}
         \le C(1+t^{-1/2})e^{-\gamma t}.
    \end{equation*}
    The nonlinear remainder
    \begin{equation*}
        \mathcal N_0(Y)=\big(y^2+y^2z,-y^2-y^2z\big)^{\top},\qquad Y=(y,z)^{\top},
    \end{equation*}
    satisfies the quadratic-cubic estimates of Lemma~\ref{lem:N}.

    Consider the weighted space
    \begin{equation*}
        \mathcal X_\gamma=\left\{Y\in C([0,\infty);X):\norm{Y}_{\mathcal X_\gamma}:=\sup_{t\ge0}e^{\gamma t}\norm{Y(t)}_X<\infty\right\}.
    \end{equation*}
    For $\norm{Y}_{\mathcal X_\gamma}\le R$, we have
    \begin{equation*}
        \norm{\mathcal N_0(Y(s))}_Z \le C(R^2+R^3)e^{-2\gamma s}.
    \end{equation*}
    Therefore,
    \begin{equation*}
        e^{\gamma t}\norm{\int_0^t e^{(t-s)\mathcal A_{\mathrm{cl}}} \mathcal N_0(Y(s))\,ds}_X \le C(R^2+R^3)\int_0^t\big(1+(t-s)^{-1/2}\big)e^{-\gamma s}\,ds.
    \end{equation*}
    Moreover, for $\norm{Y_i}_{\mathcal X_\gamma}\le R$, $i=1,2$, we have
    \begin{equation*}
        \norm{\mathcal N_0(Y_1(s))-\mathcal N_0(Y_2(s))}_Z \le C(R+R^2)e^{-2\gamma s} \norm{Y_1-Y_2}_{\mathcal X_\gamma}.
    \end{equation*}
    Since
    \begin{equation*}
        \sup_{t\ge0}\int_0^t \big(1+(t-s)^{-1/2} \big)e^{-\gamma s}\,ds\le \frac{2}{\gamma}+2<\infty,
    \end{equation*}
    the contraction argument used in the proof of Theorem~\ref{thm:main} yields a unique global mild solution satisfying
    \begin{equation*}
        \norm{Y(t)}_X\le Ce^{-\gamma t}\norm{Y_0}_X,\qquad t\ge0,
    \end{equation*}
    whenever $\norm{Y_0}_X$ is sufficiently small. This proves the attainability of the decay rates in assertions (1)-(4).

    It remains to establish the maximality of the decay rates in assertions (1)-(3).

    Suppose first that the control acts only on the first equation and that the nonlinear system is locally exponentially stabilizable with some decay rate $\gamma>B$. Choose $y_0=0$ and $z_0=\varepsilon$, where $\varepsilon>0$ is sufficiently small. By the assumed exponential stabilization estimate,
    \begin{equation*}
        \norm{Y(t)}_X\le C\varepsilon e^{-\gamma t}, \qquad t\ge0,
    \end{equation*}
    where $C>0$ is independent of $\varepsilon$. Applying the spatial integration functional to the mild formulation of the second equation and using the fact that the Neumann heat semigroup preserves spatial integrals, we obtain
    \begin{equation*}
        \int_\Om z(t,x)\,dx=\varepsilon|\Om|e^{-Bt} -\int_0^t e^{-B(t-s)}\int_\Om(y^2+y^2z)(s,x)\,dx\,ds.
    \end{equation*}
    Multiplying by $e^{Bt}$ gives
    \begin{equation*}
        e^{Bt}\int_\Om z(t,x)\,dx=\varepsilon|\Om| -\int_0^t e^{Bs}\int_\Om(y^2+y^2z)(s,x)\,dx\,ds.
    \end{equation*}
    Since $\gamma>B$, the assumed exponential decay implies
    \begin{equation*}
        \lim_{t\to\infty}e^{Bt}\int_\Om z(t,x)\,dx=0.
    \end{equation*}
    Moreover, by the nonlinear estimates,
    \begin{equation*}
        \left|\int_\Om(y^2+y^2z)(s,x)\,dx\right| \le C\varepsilon^2e^{-2\gamma s}.
    \end{equation*}
    Passing to the limit as $t\to\infty$, we obtain
    \begin{equation*}
        \varepsilon|\Om|\le C\varepsilon^2 \int_0^\infty e^{-(2\gamma-B)s}\,ds =\frac{C\varepsilon^2}{2\gamma-B}.
    \end{equation*}
    This is impossible for sufficiently small $\varepsilon>0$. Hence no decay rate larger than $B$ can be achieved by a feedback control acting only on the first equation.

    Similarly, if the control acts only on the second equation, choosing $y_0=\varepsilon$, $z_0=0$, and applying the spatial integration functional to the mild formulation of the first equation, the same argument with $A$ in place of $B$ shows that no decay rate larger than $A$ can be achieved.

    In particular, taking $h_1=h_2=0$, we deduce that the decay rate of the uncontrolled system cannot exceed $\min\{A, B\}$. This proves the maximality of the decay rates in assertions (1)-(3).

    Consequently, local exponential stabilization of the trivial steady state at a decay rate $\gamma>\max\{A,B\}$ requires controls acting on both equations. This completes the proof.
\end{proof}
\begin{remark}[Number of actuators]
    The phrase ``control acting on one equation'' refers to the component in which the distributed control enters. The actuator space itself may have dimension $m>1$. In general, the required dimension is constrained by geometric multiplicities of critical eigenvalues.
\end{remark}

\begin{remark}[Constant and variable coefficients]
    For a nonzero constant equilibrium, all coefficients of $\cG$ are constant and one may diagonalize the diffusion part by Neumann eigenfunctions. For a spatially varying equilibrium, the multiplication operators by $\overline u^2$ and $\overline u \, \overline v$ do not commute with the Neumann Laplacian; the full operator $\cG$ is generally non-self-adjoint. The Riesz projection and the Fattorini-Hautus argument above replace the scalar modal decomposition.
\end{remark}
\section{Conclusion}
\label{sec:conclusion}
    We established local exponential stabilization of every prescribed steady state of the Gray-Scott system, other than $(0,1)^\top$, at any prescribed decay rate by finite-dimensional internal feedback controls supported in an arbitrary nonempty open subset of the domain and acting only on the first equation. The same result holds when the control acts only on the second equation. 

    For a spatially varying steady state, the linearized operator has variable coefficients and is generally non-self-adjoint. We therefore isolate its finite-dimensional critical spectral subspace by means of a Riesz projection. The nonvanishing property of the steady-state coupling, together with Aronszajn's unique continuation theorem, yields the Fattorini-Hautus condition. Controllability of the critical finite-dimensional dynamics enables the construction of a Riccati-based feedback, while analytic-semigroup smoothing and the quadratic-cubic structure of the nonlinear remainder yield local exponential stabilization of the nonlinear system in $H^1(\Om)^2$.

    The trivial steady state $(0,1)^\top$ is locally exponentially stable without control, with maximal decay rate $\min\{A,B\}$. When the control acts only on the first or second equation, the maximal achievable decay rates are $B$ and $A$, respectively. Thus, in contrast to the nontrivial steady states, local exponential stabilization of the trivial steady state at an arbitrary prescribed decay rate requires controls acting on both equations.

    A natural direction for future research is to investigate whether analogous stabilization results can be established by finite-dimensional boundary feedback controls acting on only one component of the Gray-Scott system.

\section*{Acknowledgements}
    The authors thank the Department of Mathematics at the Indian Institute of Technology Kanpur for providing a conducive research environment. A.~Patel acknowledges financial support from the Ministry of Education, Government of India, through the institute fellowship. 


\end{document}